\documentclass[12pt]{article}
\usepackage{amsxtra,amssymb,amsthm,amsmath,latexsym}

\theoremstyle{plain}

\def\oH{\buildrel\circ\over H}
\def\oH1{\buildrel\circ\over H\kern-.02in{}^1}

\begin{document}


\title{ Necessary and sufficient condition for compactness of the
embedding operator
   \thanks{key words: Banach spaces, compactness, embedding operator
    }
   \thanks{AMS subject classification: 46B50, 46E30, 47B07}
}

\author{
A.G. Ramm\\
Mathematics Department, 
Kansas State University, \\
 Manhattan, KS 66506-2602, USA\\
ramm@math.ksu.edu\\
}

\date{}

\maketitle\thispagestyle{empty}

\begin{abstract} 
An improvement of the author's 
result, proved in 1961, concerning necessary and sufficient 
conditions for
the compactness of an imbedding operator is given.  
\end{abstract}


\section{Introduction}

The basic result of this note is:

{\bf Theorem 1.} {\it Let $X_1\subset X_2\subset X_3$ be Banach spaces,
$||u||_1\geq ||u||_2\geq ||u||_3$ (i.e., the norms are comparable) and  
if $||u_n||_3\to 0$ as $n\to \infty$ and $u_n$ is fundamental in $X_2$, 
 then 
$||u_n||_2\to 0$,  (i.e., the norms in $X_2$ and $X_3$ are compatible).
Under the above assumptions the embedding operator $i: X_1\to X_2$ is 
compact if and only if
the following two conditions are valid:

a) The embedding operator  $j: X_1\to X_3$ is compact,

and the following inequality holds:

b) $||u||_2\leq s ||u||_1 + c(s)||u||_3, \,\,\,\forall  
u\in X_1,$ $\forall s\in (0,1)$, where $c(s)>0$ is a constant.}

This result is an improvement of the author's old result,
proved in 1961 (see [1]), where $X_2$ 
was assumed to be a Hilbert space. The proof of Theorem 1 is simpler than 
the one in [1]. 

\section{Proof}

1. Assume that a) and b) hold and let us prove the compactness of 
$i$. Let $S=\{u: u\in X_1, ||u||_1=1\}$ be the unit sphere in $X_1$.
Using assumption a), select a sequence $u_n$ which converges 
in $ X_3$. 
We claim that this sequence 
converges also in $ X_2$. Indeed, since $||u_n||_1=1$, one 
uses assumption b) to get 
$$||u_n-u_m||_2\leq
s||u_n-u_m||_1+c(s)||u_n-u_m||_3\leq 2s +c(s)||u_n-u_m||_3.$$
Let $\eta>0$ be an arbitrary small given number. Choose $s>0$ such that 
$2s<\frac 1 2\eta$, and for a fixed $s$ choose $n$ and $m$ so large that
$c(s)||u_n-u_m||_3<\frac 1 2\eta$. This is possible because the
sequence $u_n$ converges in $ X_3$. Consequently, 
 $||u_n-u_m||_2\leq \eta$ if  $n$ and $m$ are sufficiently large.
This means that the sequence $u_n$  converges in $ X_2$.
Thus, the embedding  $i: X_1\to X_2$ is compact.
In the above argument the compatibility of the norms was not 
used. 

2. Assume now that  $i$ is compact. Let us prove that 
assumptions a) and b) hold.
Assumption a) holds because  $||u||_2\geq ||u||_3$. 
Suppose that assumption b) fails.
Then there is a sequence $u_n$ and a number $s_0>0$ such 
that $||u_n||_1=1$ and
$$
||u_n||_2\geq s_0+n||u_n||_3.
\eqno{(1)}$$
If  the embedding operator $i$ is compact and $||u_n||_1=1$, 
then one may assume that
the sequence $u_n$ converges in $X_2$. Its limit cannot be equal to zero, 
because, by (1), $||u_n||_2\geq s_0>0$. The  sequence $u_n$ converges in 
$X_3$  because $||u_n-u_m||_2\geq ||u_n-u_m||_3$, and its 
limit in $X_3$ is not
zero, because the norms in $X_3$ and in $X_2$ are compatible.
Thus, (1) implies $||u_n||_1=O(\frac 1 n)\to 0$ as $n\to \infty$,
while $||u_n||_1=1$. This is a contradiction, which proves 
that b) holds.

Theorem 1 is proved. \hfill $\Box$

\end{document}